\documentclass{siamltex}

\usepackage{amsfonts}
\usepackage{amsmath}
\usepackage{mathrsfs}

\newcommand{\bdW}[1]{\overleftarrow{dW}_{\!\! #1}}

\title{Almost Global Stochastic Stability\thanks{This work was supported 
by the ARO under Grant DAAD19-03-1-0073.}}

\author{Ramon van Handel\thanks{The author is with the departments of 
Physics and Control \& Dynamical Systems, California Institute of 
Technology 266-33, Pasadena, CA 91125, USA (ramon@its.caltech.edu).}}

\begin{document}
\maketitle

\begin{abstract}
	We develop a method to prove almost global stability of stochastic
	differential equations in the sense that almost every initial 
	point (with respect to the Lebesgue measure) is asymptotically
	attracted to the origin with unit probability.  The method can be 
	viewed as a dual to Lyapunov's second method for stochastic 
	differential equations and extends the deterministic result in 
	[A.\ Rantzer, {\em Syst. Contr. Lett.}, 42 (2001), pp.\ 161--168].
	The result can also be used in certain cases to find stabilizing
	controllers for stochastic nonlinear systems using convex 
	optimization.  The main technical tool is the theory of stochastic 
	flows of diffeomorphisms.
\end{abstract}

\begin{keywords} 
	stochastic stability,
	stochastic flows,
	nonlinear stochastic control
\end{keywords}

\begin{AMS}
34F05, 60H10, 93C10, 93D15, 93E15
\end{AMS}

\pagestyle{myheadings}
\thispagestyle{plain}
\markboth{RAMON VAN HANDEL}{ALMOST GLOBAL STOCHASTIC STABILITY}

\section{Introduction}

Lyapunov's second method or the method of Lyapunov functions, though 
developed in the late 19th century, remains one of the most important 
tools in the study of deterministic differential equations.
The power of the method lies in the fact that an important qualitative 
property of a differential equation, the stability of an equilibrium 
point, can be proved without solving the equation explicitly. The theory 
was generalized to stochastic differential equations in the 1960s with 
fundamental contributions by Has'minski\u{\i} \cite{hasminskii} and 
Kushner \cite{kushner1}.

Lyapunov's method also underlies many important applications in the area 
of nonlinear control \cite{isidori}.  Finding optimal controls for 
nonlinear systems is generally an intractable problem, but often a 
solution can be found which stabilizes the system.  Unlike in 
deterministic control theory, where nonlinear control is now a major 
field, there are very few results on stochastic nonlinear control.  It is 
only recently that stochastic versions of the classical stabilization 
results of Jurdjevic-Quinn, Artstein and Sontag were developed by 
Florchinger \cite{florch1,florch2,florch4} and backstepping designs for 
stochastic strict-feedback systems were developed by Deng and Krsti\'c 
\cite{deng1,deng2}.

In this paper we will not consider stochastic stability in the sense of 
Has'minski\u{\i}; rather, we ask the following question: for a given 
It\^o stochastic differential equation on $\mathbb{R}^n$, can we prove 
that for almost every (with respect to the Lebesgue measure on 
$\mathbb{R}^n$) initial state the solution of the equation converges to 
the origin almost surely as $t\to\infty$, i.e., is the origin {\em almost} 
globally stable?  This notion of stability is clearly weaker than global 
stability in the sense of Has'minski\u{\i}, but is of potential interest 
in many cases in which global stability may not be attained.

Our main result is a Lyapunov-type theorem that can be used to prove
almost global stability of stochastic differential equations, extending
the deterministic result of Rantzer \cite{rantzer}.  The theorem has
several remarkable properties.  It can be viewed as a ``dual'' to
Lyapunov's second method in the following sense: whereas the Lyapunov
condition reads $\mathscr{L}V<0$, where $\mathscr{L}$ is the
characteristic operator of the stochastic differential equation and $V$ is
the Lyapunov function, the condition that guarantees almost global
stability reads $\mathscr{L}^*D<0$ where $\mathscr{L}^*$ is the formal
adjoint of $\mathscr{L}$ (also known as the Fokker-Planck operator.)  
Hence the relation between the two theorems recalls the duality between
densities and expectations which is prevalent throughout the theory of
stochastic processes.

A further interesting property is the following convexity property.
Suppose we are given an It\^o equation of the form
$$
	x_t = x+\int_s^t(X_0(x_\tau)+u(x_\tau)Y(x_\tau))\,d\tau+
		\sum_{k=1}^m\int_s^t X_k(x_\tau)\,dW_\tau^k
$$
where $X_k(0)=0$ for $k=0\ldots m$ and $u(x)$ is a state feedback control.
The goal is to design $u(x)$ such that the origin is almost globally 
stable.  It is easily verified that the set of pairs of functions 
$(D(x),u(x)D(x))$ which satisfy $\mathscr{L}^*D<0$ is convex.  Note that 
the classical Lyapunov condition $\mathscr{L}V<0$ is not convex.

The above convexity property was used in the deterministic case by Prajna, 
Parrilo and Rantzer \cite{prajna1} to formulate the search for almost 
globally stabilizing controllers as a convex optimization problem, 
provided that $X_k$, $Y$, $D$ and $u$ are rational functions.  The method 
applies equally to the stochastic case and thus provides a tool for 
computer-aided design of stochastic nonlinear controllers.

It must be emphasized that almost global stability is a global property of 
the flow which places very few restrictions on the local behavior near the 
origin.  In particular, local stability is not implied\footnote{
	The term ``stability'' seems a bit of a misnomer; despite that 
	almost all points converge to the origin, a trajectory that 
	starts close to the origin could move very far from the origin
	before converging to it.  We have used the term that has been
	used in the deterministic literature, e.g.\ \cite{manzon}.
}.  A very fruitful approach to studying the local dynamical behavior of 
stochastic differential equations (and more general random dynamical 
systems) is developed by Arnold \cite{arnold3}.  First, the flow 
associated to the stochastic differential equation is linearized; then 
Oseledec's multiplicative ergodic theorem is used to provide a suitable 
``time-averaged'' notion of the eigenvalues of the linearized flow.  To 
prove almost global stability we do not linearize the flow, though the 
proofs still rely on the flow of diffeomorphisms generated by the 
stochastic equation.  We refer to \cite{arnold2} for an introduction to 
the dynamical approach to stochastic analysis.

This paper is organized as follows.  In \S\ref{sec:not} we fix the
notation that will be used in the remainder of the paper.  In
\S\ref{sec:det} we reproduce the deterministic result of Rantzer
\cite{rantzer} with a significantly different proof that generalizes to
the stochastic case.  In \S\ref{sec:flow} we review the theory of
stochastic flows of diffeomorphisms generated by stochastic differential
equations.  \S\ref{sec:main} is devoted to the statement and proof of our
main result for the case of globally Lipschitz continuous coefficients.  
In \S\ref{sec:further} the main result is extended to cases in which the
global Lipschitz condition does not neccessarily hold. A few examples are
given in \S\ref{sec:examples}.  Finally, in \S\ref{sec:control} we discuss
the application to control synthesis using convex optimization.

\section{Notation}
\label{sec:not}

Throughout this article we will consider (stochastic) differential 
equations in $\mathbb{R}^n$.  The Lebesgue measure on $\mathbb{R}^n$ will 
be denoted by $\mu$.  $\mathbb{R}_+$ denotes the nonnegative real numbers 
and $\mathbb{Z}_+$ the nonnegative integers.

We remind the reader of the following definitions: for $0<\alpha\le 1$, a 
function $f:X\to Y$ from a normed space $(X,\|\cdot\|)$ to a normed 
space $(Y,\|\cdot\|_Y)$ is called globally H\"older continuous of order 
$\alpha$ if there exists a positive constant $C$ such that
\begin{equation}
\label{eq:holder}
	\|f(x)-f(y)\|_Y\le C\|x-y\|^\alpha~~~~~~~\forall x,y\in X.
\end{equation}
$f$ is locally H\"older continuous of order $\alpha$ if it satisfies the 
condition (\ref{eq:holder}) on every bounded subset of $X$.  $f$ is called 
globally (locally) Lipschitz continuous if it is globally (locally) 
H\"older continuous of order $1$.  $f$ is called a $C^{k,\alpha}$ function 
if it is $k$ times continuously differentiable and the $k$-th derivatives 
are locally H\"older continuous of order $\alpha$ for some 
$k\in\mathbb{Z}_+$ and $0<\alpha\le 1$.

\section{The deterministic case}
\label{sec:det}

In this section we give a new proof of Rantzer's theorem \cite{rantzer}
which is a deterministic counterpart of our main result.  Our proof 
demonstrates the main features of the proof of the stochastic result in 
the simpler deterministic case.

The following Lemma is similar to Lemma A.1 in \cite{rantzer}, and we omit 
the proof.

\begin{lemma}
\label{lemma:detliou}
	Let $f\in C^1(\mathbb{R}^n,\mathbb{R}^n)$ be globally Lipschitz 
	continuous, $S\subset\mathbb{R}^n$ be an invariant set of $\dot 
	x(t)=f(x(t))$ and $Z\subset S$ be $\mu$-measurable. Let $D\in 
	C^1(S,\mathbb{R})$ be integrable on $Z$.  Then
	\begin{equation}
	\label{eq:detliou}
		\int_{\phi_t^{-1}(Z)}D(x)\,dx=
		\int_ZD(x)\,dx-\int_0^t\int_{\phi_\tau^{-1}(Z)}
			[\nabla\cdot(fD)](x)\,dx\,d\tau
	\end{equation}
	where $\phi_t:\mathbb{R}^n\to\mathbb{R}^n$ is the flow
	of $f$.
\end{lemma}

\begin{theorem}
\label{thm:determin}
	Let $f\in C^1(\mathbb{R}^n,\mathbb{R}^n)$ be globally Lipschitz 
	continuous and let $f(0)=0$.  Suppose there exists
	$D\in C^1(\mathbb{R}^n\backslash\{0\},\mathbb{R}_+)$ 
	such that $D$ is integrable on $\{x\in\mathbb{R}^n:|x|>1\}$ and
	$[\nabla\cdot(fD)](x)>0$ for $\mu$-almost all $x$.  Then for 
	$\mu$-almost all initial states $x(0)$ the solution of $\dot 
	x(t)=f(x(t))$ tends to the origin as $t\to\infty$.
\end{theorem}

\begin{proof}
Let $S=\mathbb{R}^n\backslash\{0\}$, $\varepsilon>0$ and 
$Z=\{x\in\mathbb{R}^n:|x|>\varepsilon\}$.  Note that $\phi_t(x)$ is a 
diffeomorphism for every $t\in\mathbb{R}$; hence $\phi_t(x)$ is 
one-to-one, and as $\phi_t(0)=0$, $t\in\mathbb{R}$ is a solution of $\dot 
x(t)=f(x(t))$ there can be no $x\in S$ such that $\phi_t(x)=0$ for some 
$t\in\mathbb{R}$.  We have thus verified the invariance of $S$ under the 
flow $\phi_t(x)$.  We now invoke Lemma \ref{lemma:detliou}.  As $D(x)$ is 
nonnegative expression (\ref{eq:detliou}) is also nonnegative.  
Furthermore, (\ref{eq:detliou}) is finite because $D$ is integrable on 
$Z$, and is nonincreasing due to $[\nabla\cdot(fD)](x)\ge 0$.  By monotone 
convergence the limit as $t\to\infty$ exists and is finite.  Hence
$$
	\int_0^\infty\mathbb{D}(\phi_\tau^{-1}(Z))\,d\tau
	<\infty,
	~~~~~~~ ~~~~~~~
	\mathbb{D}(A)=
	\int_{A}[\nabla\cdot(fD)](x)\,dx.
$$
Note that the assumption $[\nabla\cdot(fD)](x)\ge 0$ implies that 
$\mathbb{D}$ is a measure on $S$.  The measure space 
$(S,\mathbb{D})$ is $\sigma$-finite as 
$\mathbb{D}(\{x\in S:\tfrac{1}{k}<|x|<k\})<\infty$ for all $k>1$ and 
$\bigcup_{n=2}^\infty\{x\in S:\tfrac{1}{k}<|x|<k\}=S$.

We now fix some $m\in\mathbb{N}$ and divide the halfline into bins 
$S_k^m=[(k-1)2^{-m},k2^{-m}]$, $k\in\mathbb{N}$.  From each bin we choose 
a time $t_k^m\in S_k^m$ such that
$$
	\mathbb{D}(\phi_{t_k^m}^{-1}(Z))\le
		\inf_{t\in S_k^m}
		\mathbb{D}(\phi_{t}^{-1}(Z))
		+2^{-k}.
$$
For fixed $m$, we denote this discrete grid by 
$T_m=\{t_k^m:k\in\mathbb{N}\}$.  We now have
$$
	2^{-m}\sum_{k=1}^\infty \mathbb{D}(\phi_{t_k^m}^{-1}(Z))
	\le 2^{-m}+
	\int_0^\infty\mathbb{D}(\phi_\tau^{-1}(Z))\,d\tau
	<\infty.
$$
As $\mathbb{D}$ is $\sigma$-finite we can now apply the Borel-Cantelli 
lemma, which gives
$$
	\mathbb{D}\left(\limsup_{k\to\infty}\phi_{t_k^m}^{-1}(Z)\right)
	=\mu\left(\limsup_{k\to\infty}\phi_{t_k^m}^{-1}(Z)\right)=0
$$
where the first equality follows as $[\nabla\cdot(fD)](x)>0$ 
$\mu$-a.e.\ implies $\mu\ll\mathbb{D}$.  Consequently
\begin{equation*}
        \mu\left(\bigcup_{m=1}^\infty\limsup_{t\in T_m}\phi_{t}^{-1}(Z)\right)
        \le \sum_{m=1}^\infty
        \mu\left(\limsup_{t\in T_m}\phi_{t}^{-1}(Z)\right)=0.
\end{equation*}
We have thus shown that the set of initial states $x$ for which there 
are, for some $m$, infinitely many times $t\in T_m$ such that 
$\phi_t(x)\in Z$, has Lebesgue measure zero.

We now claim that if $\limsup_{t\to\infty}|\phi_t(x)|>\varepsilon$, then 
we can choose $m$ so that there are infinitely many times $t$ in $T_m$ 
such that $\phi_t(x)\in Z$.  The statement is trivial if also
$\liminf_{t\to\infty}|\phi_t(x)|>\varepsilon$; let us thus assume 
that $\liminf_{t\to\infty}|\phi_t(x)|\le\varepsilon$.  We will need the 
following result.  Due to the global Lipschitz condition and $f(0)=0$, 
we have
\begin{equation}
\label{eq:determinarg}
	|\phi_t(x)|\le |\phi_s(x)|+\int_s^t|f(\phi_\sigma(x))|\,d\sigma
		\le |\phi_s(x)|+C\int_s^t|\phi_\sigma(x)|\,d\sigma
\end{equation}
for some constant $C>0$.  Thus Gronwall's lemma gives
$|\phi_t(x)|\le |\phi_s(x)|\,e^{C(t-s)}$.  Now note that
$\liminf_{t\to\infty}|\phi_t(x)|\le\varepsilon<\limsup_{t\to\infty}|\phi_t(x)|$
implies that there exist $\varepsilon''>\varepsilon'>\varepsilon$ such 
that (i) there are infinitely many upcrossings of the curve 
$|\phi_t(x)|$ through $\varepsilon'$, and (ii) $|\phi_t(x)|$ crosses 
$\varepsilon''$ infinitely often.  Denote by $t''$ a time such that 
$|\phi_{t''}(x)|=\varepsilon''$ and by $t'$ the latest time previous to 
$t''$ that $|\phi_{t'}(x)|=\varepsilon'$.  Then clearly 
$t''-t'\ge\tfrac{1}{C}\log\tfrac{\varepsilon''}{\varepsilon'}$.  As this 
happens infinitely often, we conclude that $\phi_{t}(x)$ infinitely 
often spends a time in excess of 
$\tfrac{1}{C}\log\tfrac{\varepsilon''}{\varepsilon'}$ in $Z$.  But then 
clearly $m$ can be chosen large enough so that every such interval 
includes at least one of the $t_k^m\in T_m$.

We have now shown that for $\mu$-almost all $x\in\mathbb{R}^n$ we have 
$\limsup_{t\to\infty}|\phi_t(x)|\le\varepsilon$, i.e.\ for $\mu$-almost 
all $x\in\mathbb{R}^n$ $\exists t_e>0$ such that $|\phi_t(x)|\le\varepsilon$
for $t\ge t_e$.  But as this holds for any $\varepsilon>0$ the 
trajectories must converge to the origin.  \qquad
\end{proof}

\section{Stochastic flows}
\label{sec:flow}

The purpose of this section is to review, without proofs, some results of 
the theory of stochastic flows of diffeomorphisms generated by stochastic 
differential equations.  A detailed exposition on the subject can be found 
in \cite{kunita1,kunita2} and shorter treatments are in
\cite{arnold2,arnold3,baxendale,kunita3}.

Throughout this article $(\Omega,\mathcal{F},\mathbb{P})$ denotes the
canonical Wiener space of the $m$-dimensional Brownian motion $W_t$ with
two-sided time $\mathbb{R}$.  We also introduce the two-parameter
filtration $\mathcal{F}_s^t=\sigma\{W^k_u-W^k_v:s\le v\le u\le t,
~1\le k\le m\}$.  The extension to two-sided time is important in that it 
allows us to treat the Wiener process as a dynamical system 
\cite{arnold2,arnold3,baxendale}:

\begin{theorem}
\label{thm:theta}
	There exists a one-parameter group $\{\theta_t:t\in\mathbb{R}\}$
	of measure-preserving transformations of 
	$(\Omega,\mathcal{F},\mathbb{P})$ such that
	$W_t(\theta_s\omega)=W_{t+s}(\omega)-W_s(\omega)$ for all
	$\omega\in\Omega$ and $s,t\in\mathbb{R}$.
\end{theorem}

We will consider It\^o stochastic differential equations of the form
\begin{equation}
\label{eq:SDE}
	x_t = x+\int_s^tX_0(x_\tau)\,d\tau+
		\sum_{k=1}^m\int_s^t X_k(x_\tau)\,dW_\tau^k
\end{equation}
with the following assumptions:
\begin{enumerate}
\item $x\in\mathbb{R}^n$.
\item $X_k:\mathbb{R}^n\to\mathbb{R}^n$, $k=0\ldots m$ are globally
	Lipschitz continuous.
\end{enumerate}
The global Lipschitz condition guarantees many nice properties of the 
solutions; we will assume it for the time being, and later relax this 
requirement somewhat in \S\ref{sec:further}.

Denote by $\xi_{s,t}(x,\omega)$ (or simply $\xi_{s,t}(x)$) the solution of
(\ref{eq:SDE}) at time $t\ge s$ given the initial condition $x_s=x$.  It 
is well known that in the case of globally Lipschitz continuous 
coefficients there exists a unique, non-exploding solution $\xi_{s,t}(x)$ 
which is an $\mathcal{F}_s^t$-semimartingale and is in $L^p$ for any 
$p\ge 1$ (e.g.\ \cite{kunita1}).

\begin{theorem}[\cite{kunita1,arnold3}]
	Suppose $X_k$, $k=0\ldots m$ are globally Lipschitz
	continuous and let $s<t$.  Then we have the following 
	properties:
	\begin{enumerate}
	\item	$\xi_{s,s}(x,\omega)=x$ for all $s$ and $\omega$.
	\item	For any $u$ we have $\xi_{s,t}(\cdot,\theta_u\omega)=
		\xi_{s+u,t+u}(\cdot,\omega)$.
	\item	For almost all $\omega$ we have
		$\xi_{s,t}(\cdot,\omega)=
		\xi_{r,t}(\xi_{s,r}(\cdot,\omega),\omega)$ for all $s<r<t$.
	\item	$\xi_{s,t}(x)$ is $\mathbb{P}$-a.s.\ continuous in
		$(s,t,x)$.
	\item   For a.e.\ $\omega$ the map
		$\xi_{s,t}(\cdot,\omega):\mathbb{R}^n\to\mathbb{R}^n$
		is a homeomorphism for all $s<t$.
	\end{enumerate}
\end{theorem}

The following result establishes that, under additional smoothness
conditions, $\xi_{s,t}(x)$ is in fact a {\em stochastic flow of
diffeomorphisms}.

\begin{theorem}[\cite{kunita1}]
	Suppose $X_k$, $k=0\ldots m$ are globally Lipschitz continuous
	and that they are $C^{p,\alpha}$ functions for some $p\ge 1$ and 
	$0<\alpha<1$.  Then for almost all $\omega$ the map 
	$\xi_{s,t}(\cdot,\omega):\mathbb{R}^n\to\mathbb{R}^n$ is a $C^p$ 
	diffeomorphism for any $s\le t$, and
	\begin{multline}
	\label{eq:diffflowf}
		\frac{\partial\xi_{s,t}(x)^i}{\partial x^j}=
		\delta^i_j+\sum_{\beta=1}^n\int_s^t
		\frac{\partial X_0^i}{\partial x^\beta}
		(\xi_{s,\tau}(x))\,
		\frac{\partial\xi_{s,\tau}(x)^\beta}{\partial x^j}
		\,d\tau \\
		+\sum_{k=1}^m\sum_{\beta=1}^n\int_s^t
		\frac{\partial X_k^i}{\partial x^\beta}
		(\xi_{s,\tau}(x))\,
		\frac{\partial\xi_{s,\tau}(x)^\beta}{\partial x^j}
		\,dW_\tau^k.
	\end{multline}
\end{theorem}

It will be convenient for our purposes to work with the inverse flow 
$\xi_{s,t}^{-1}(x)$, considered as a backward stochastic process in the 
time variable $s$ (with $t$ fixed).  This will not give rise to ordinary 
It\^o integrals as $s$ behaves like a time-reversed variable, and hence 
the adaptedness of the process runs backwards in time.  The It\^o 
backward integral is defined as \cite{kunita1}
\begin{equation*}
	\int_s^t f_\sigma\bdW{\sigma}\equiv
	\mbox{lim in prob}\,\sum_{k=0}^{n-1}f_{t_{k+1}}(W_{t_{k+1}}-
		W_{t_k})
\end{equation*}
where $f_s$ is a backward predictable process with
$\int_s^t|f_u|^2\,du<\infty$ a.s., and the formal construction of the
integral from simple functions proceeds along the usual lines.  The
backward integral has similar properties to the forward integral; in
particular, it is a backward $\mathcal{F}_s^t$-local martingale (for fixed 
$t$) and satisfies an It\^o formula (e.g.\ \cite{bensoussan}, pp.\ 124) 
which is proved in the same way as its forward counterpart: given
$\xi_s=\xi_t+\int_s^ta_\sigma\,d\sigma+\sum_k\int_s^t(b_\sigma)_k\,
\bdW{\sigma}^k$ with backward predictable processes $a_s$, $(b_s)_k$ s.t.\
$\int_s^ta_\sigma\,d\sigma<\infty$ a.s., 
$\int_s^t|(b_\sigma)_k|^2\,d\sigma<\infty$ a.s., then for any $C^2$ 
function $F:\mathbb{R}^n\to\mathbb{R}$
\begin{multline}
\label{eq:backito}
	F(\xi_s)=F(\xi_t)+
	\frac{1}{2}\sum_{k=1}^m\sum_{i,j=1}^n\int_s^t
		(b_\sigma)_k^i(b_\sigma)_k^j
		\frac{\partial^2F}{\partial x^i\partial x^j}
			(\xi_\sigma)\,d\sigma
	\\ 
	+\sum_{i=1}^n \int_s^t
		a_\sigma^i
		\frac{\partial F}{\partial x^i}
		(\xi_\sigma)\,d\sigma
	+\sum_{k=1}^m\sum_{i=1}^n \int_s^t
		(b_\sigma)_k^i
		\frac{\partial F}{\partial 
		x^i}(\xi_\sigma)\,\bdW{\sigma}^k.
\end{multline}
We can now formulate the following result.

\begin{theorem}[\cite{kunita1}]
\label{thm:invflow}
	Suppose $X_k$, $k=0\ldots m$ are globally Lipschitz continuous
	and that they are $C^{p,\alpha}$ functions for some $p\ge 2$ and 
	$0<\alpha<1$.  Then
	\begin{equation*}
		\xi_{s,t}^{-1}(x)=x-
		\int_s^t\tilde X_0(\xi_{\sigma,t}^{-1}(x))\,d\sigma-
		\sum_{k=1}^m\int_s^t
		X_k(\xi_{\sigma,t}^{-1}(x))\,\bdW{\sigma}
	\end{equation*}
	where we have defined
	\begin{equation*}
		\tilde X_0(x)=X_0(x)-\sum_{k=1}^m
		\sum_{\beta=1}^nX_k^\beta(x)
			\frac{\partial}{\partial x^\beta}X_k(x).
	\end{equation*}
\end{theorem}

This expression can be manipulated much in the same way as its forward 
counterpart.  In particular, under the conditions of Theorem 
\ref{thm:invflow} and using (\ref{eq:backito}), we obtain for any $C^2$ 
function $F:\mathbb{R}^n\to\mathbb{R}$ the backward It\^o formula
\begin{multline}
\label{eq:fbackito}
	F(\xi_{s,t}^{-1}(x))=F(x)+
	\frac{1}{2}\sum_{k,i,j}\int_s^t
		X_k^i(\xi_{\sigma,t}^{-1}(x))
		X_k^j(\xi_{\sigma,t}^{-1}(x))
		\frac{\partial^2F}{\partial x^i\partial x^j}
			(\xi_{\sigma,t}^{-1}(x))\,d\sigma
	\\ 
	-\sum_i \int_s^t
		\tilde X_0^i(\xi_{\sigma,t}^{-1}(x))
		\frac{\partial F}{\partial x^i}
		(\xi_{\sigma,t}^{-1}(x))\,d\sigma
	-\sum_{k,i} \int_s^t
		X_k^i(\xi_{\sigma,t}^{-1}(x))
		\frac{\partial F}{\partial 
		x^i}(\xi_{\sigma,t}^{-1}(x))\,\bdW{\sigma}^k.
\end{multline}
Similarly we can differentiate the inverse flow, giving
\begin{multline}
\label{eq:backjaco}
	\frac{\partial\xi_{s,t}^{-1}(x)^i}{\partial x^j}=
	\delta^i_j-\sum_{\beta=1}^n\int_s^t
	\frac{\partial\tilde X_0^i}{\partial x^\beta}
	(\xi_{\sigma,t}^{-1}(x))\,
	\frac{\partial\xi_{\sigma,t}^{-1}(x)^\beta}{\partial x^j}
	\,d\sigma \\
	-\sum_{k=1}^m\sum_{\beta=1}^n\int_s^t
	\frac{\partial X_k^i}{\partial x^\beta}
	(\xi_{\sigma,t}^{-1}(x))\,
	\frac{\partial\xi_{\sigma,t}^{-1}(x)^\beta}{\partial x^j}
	\,\bdW{\sigma}^k.
\end{multline}
This expression is obtained, in the same way as its forward counterpart 
(\ref{eq:diffflowf}), by letting $y\to 0$ in the backward expression 
corresponding to \cite{kunita1}, pp.\ 219, Eq.\ (4).

\section{The main result}
\label{sec:main}

We consider an It\^o equation of the form (\ref{eq:SDE}).
We write
\begin{equation*}
	\mathscr{L}^*f(x)=
	\frac{1}{2}\sum_{k=1}^m\sum_{i,j=1}^n
		\frac{\partial^2}{\partial x^i\partial x^j}
		(X^i_k(x)X^j_k(x)f(x))
	-\sum_{i=1}^n
		\frac{\partial}{\partial x^i}
		(X_0^i(x)f(x)).
\end{equation*}
The following is our main result:

\begin{theorem}
\label{thm:main}
Let $X_k$, $k=0\ldots m$ be globally Lipschitz continuous and 
$C^{2,\alpha}$ for some $\alpha>0$, and let $X_k(0)=0$.
Suppose there exists $D\in C^2(\mathbb{R}^n\backslash\{0\},\mathbb{R}_+)$
such that $D$ is integrable on $\{x\in\mathbb{R}^n:|x|>1\}$ and
$\mathscr{L}^*D(x)<0$ for $\mu$-almost all $x$.  Then for every initial 
time $s$ and $\mu$-almost every initial state $x$ the flow $\xi_{s,t}(x)$ 
tends to the origin as $t\to\infty$ $\mathbb{P}$-a.s.
\end{theorem}

Before we prove the theorem, let us prove a stochastic version of Lemma 
\ref{lemma:detliou}.

\begin{lemma}
\label{lemma:stochliou}
Suppose $X_k$, $k=0\ldots m$ are globally Lipschitz continuous and 
$C^{2,\alpha}$ functions for some $\alpha>0$.  Let 
$S_\ell\subset S_{\ell+1}\subset\mathbb{R}^n$ be an increasing sequence of 
open sets such that $\tau_\ell=\sup\{s<t:\xi_{s,t}^{-1}(x)\not\in 
S_\ell\}\to-\infty$ as $\ell\to\infty$ $\mathbb{P}$-a.s.\ for every $x\in 
S=\bigcup_\ell S_\ell$.  Suppose there is a $D\in C^2(S,\mathbb{R}_+)$ 
that is integrable on a measurable set $Z\subset S$, that obeys 
$\mathscr{L}^*D\le 0$ on $S$, and such that for each $\ell$ there is a 
$D_\ell\in C^2(\mathbb{R}^n,\mathbb{R}_+)$ that coincides with $D$ on 
$S_\ell$.  Then
$$
	0\le \int_ZD(x)\,dx+
	\int_s^t\mathbb{E}
	\int_{\xi_{\sigma,t}^{-1}(Z)}
		\mathscr{L}^*D(x)\,dx\,d\sigma
$$
for all $s\le t$, and in particular the limit as $s\to-\infty$ of this
expression is well defined.
\end{lemma}

\begin{proof}
Denote by ${\bf J}_{s,t}(x)$ the matrix with elements
${\bf J}_{s,t}(x)^i_j=\partial\xi_{s,t}^{-1}(x)^i/\partial x^j$, i.e.\
$$
	{\bf J}_{s,t}(x)^i_j=\delta^i_j
	-\sum_\alpha\int_s^t
	\frac{\partial\tilde X_0^i}{\partial x^\alpha}
		(\xi_{\sigma,t}^{-1}(x))
	{\bf J}_{\sigma,t}(x)^\alpha_j\,d\sigma
	-\sum_{k,\alpha}\int_s^t
	\frac{\partial X_k^i}{\partial x^\alpha}(\xi_{\sigma,t}^{-1}(x))
	{\bf J}_{\sigma,t}(x)^\alpha_j\,\bdW{\sigma}^k
$$
by Eq.\ (\ref{eq:backjaco}).  Denote by $|{\bf J}_{s,t}(x)|$ its 
determinant, i.e.\
$$
	|{\bf J}_{s,t}(x)|=\sum_{j_1\cdots j_n=1}^n
	\varepsilon_{j_1,\ldots,j_n}
	{\bf J}_{s,t}(x)_{j_1}^1
	{\bf J}_{s,t}(x)_{j_2}^2\cdots
	{\bf J}_{s,t}(x)_{j_n}^n
$$
where $\varepsilon_{j_1,\ldots,j_n}$ is the antisymmetric tensor.  
Using It\^o's rule and straightforward calculations we obtain
\begin{equation*}
\begin{split}
	|&{\bf J}_{s,t}(x)|=1
	-\sum_i\int_s^t
		\frac{\partial\tilde X_0^i}{\partial x^i}
		(\xi_{\sigma,t}^{-1}(x))
		|{\bf J}_{\sigma,t}(x)|\,d\sigma
	-\sum_{k,i}\int_s^t
		\frac{\partial X_k^i}{\partial x^i}
		(\xi_{\sigma,t}^{-1}(x))
		|{\bf J}_{\sigma,t}(x)|\,\bdW{\sigma}^k
	 \\ &+
	\frac{1}{2}\sum_{k,i,j}\int_s^t
	\left[
	\frac{\partial X_k^i}{\partial x^i}(\xi_{\sigma,t}^{-1}(x))
	\frac{\partial X_k^j}{\partial x^j}(\xi_{\sigma,t}^{-1}(x))-
	\frac{\partial X_k^i}{\partial x^j}(\xi_{\sigma,t}^{-1}(x))
	\frac{\partial X_k^j}{\partial x^i}(\xi_{\sigma,t}^{-1}(x))
	\right]|{\bf J}_{\sigma,t}(x)|\,d\sigma.
\end{split}
\end{equation*}
Note that as $\xi_{s,t}^{-1}(\cdot)$ is a diffeomorphism a.s., its 
Jacobian ${\bf J}_{s,t}(\cdot)$ must a.s.\ be an invertible matrix; but as 
$|{\bf J}_{s,t}(x)|$ has a.s.\ continuous sample paths and $|{\bf 
J}_{t,t}(x)|=1$, this implies that a.s.\ $|{\bf J}_{s,t}(x)|>0$ for all 
$s<t$.  Using (\ref{eq:fbackito}) with $F=D_\ell$ and It\^o's rule we 
obtain
\begin{multline*}
	0\le D_\ell(\xi_{s,t}^{-1}(x))|{\bf J}_{s,t}(x)|=D_\ell(x)+
	\int_s^t(\mathscr{L}^*D_\ell)(\xi_{\sigma,t}^{-1}(x))
		|{\bf J}_{\sigma,t}(x)|\,d\sigma \\
	-\sum_{k=1}^m\sum_{i=1}^n\int_s^t
	\frac{\partial X_k^iD_\ell}{\partial x^i}(\xi_{\sigma,t}^{-1}(x))
		|{\bf J}_{\sigma,t}(x)|\,\bdW{\sigma}^k.
\end{multline*}
Now note that as $D_\ell$ coincides with $D$ on $S_\ell$, we can identify
$(\mathscr{L}^*D_\ell)(\xi_{s\vee\tau_\ell,t}^{-1}(x))=
(\mathscr{L}^*D)(\xi_{s\vee\tau_\ell,t}^{-1}(x))$ for every $\ell$.
Moreover, as the last term in the expression above is a backward local 
martingale, there exists a sequence of stopping times 
$\tau_p'\searrow-\infty$ such that the stochastic integral stopped at 
$\tau_p'$ is a martingale.  Replacing $s$ by $s\vee\tau_\ell\vee\tau_p'$ 
in the expression above and taking the expectation gives
$$
	0\le D(x)+
	\mathbb{E}\int_{s\vee\tau_\ell\vee\tau_p'}^t
		(\mathscr{L}^*D)(\xi_{\sigma,t}^{-1}(x))
		|{\bf J}_{\sigma,t}(x)|\,d\sigma.
$$
We can now let $\ell,p\to\infty$ by monotone convergence. Integrating both 
sides gives
$$
	0\le 
	\int_Z D(x)\,dx+ \int_s^t\mathbb{E}\int_Z
		(\mathscr{L}^*D)(\xi_{\sigma,t}^{-1}(x))
		|{\bf J}_{\sigma,t}(x)|\,dx\,d\sigma
$$
where we have used Tonelli's theorem to change the order of integration.  
The result follows after a change of coordinates.
\qquad
\end{proof}

The proof of Theorem \ref{thm:main} is similar to the proof of Theorem
\ref{thm:determin}.  The stochastic version of the argument following
(\ref{eq:determinarg}), however, is a little more subtle, as we do not
have a pathwise upper bound on the rate of growth of sample paths.  On the
other hand, we can establish such a bound {\em in probability} which,
together with the strong Markov property, is sufficient for our purposes;
a similar argument was used in \cite{lasalle} to the same effect.  For
this purpose we give the following Lemma, various versions of which appear
in the literature (the result below is adapted from \cite{deng2}.)

\begin{lemma}\label{lemma:uniform}
	Let $X_k$, $k=0\ldots m$ be locally Lipschitz continuous and
	$\lambda>0$.  Then
	\begin{equation*}
		\mathbb{P}
		\left[\sup_{0\le\delta\le\Delta}|\xi_{s,s+\delta}(x)-x|\ge
			\lambda\right]\le K_1\Delta+K_2\Delta^2
	\end{equation*}
	where $K_1,K_2<\infty$ are constants that depend only on 
	$\lambda$ and $|x|$.
\end{lemma}

\begin{proof}
Let ${\bf W}_t$ be the $m$-vector with elements $W_t^k$ and ${\bf 
X}(\cdot)$ be the $n\times m$-matrix with entries $X_k^i(\cdot)$, 
$k=1\ldots m$.  For $r>0$, define $B_r(x')= 
\{x\in\mathbb{R}^n:|x-x'|<r\}$, $B_r=B_r(0)$, and
\begin{equation*}
	\rho_0(r)=\sup_{|y|<r}|X_0(y)|,\qquad
	\rho_1(r)=\sup_{|y|<r}\|{\bf X}(y)\|=
		\sup_{|y|<r}{\rm tr}[{\bf X}(y)^T{\bf X}(y)]^{1/2}.
\end{equation*}
Let $\tau_r$ be the first exit time of $\xi_{s,t}(x)$ from $B_r$.  In 
\cite{deng2}, pp.\ 1240 it was established that
\begin{equation*}
	\mathbb{E}
	\left[\sup_{0\le\delta\le\Delta}
		|\xi_{s,(s+\delta)\wedge\tau_r}(x)-x|^2
	\right]\le 2\rho_0(r)^2\Delta^2+8\rho_1(r)^2\Delta.
\end{equation*}
Hence we have by Markov's inequality
\begin{equation*}
	\mathbb{P}
	\left[\sup_{0\le\delta\le\Delta}
		|\xi_{s,(s+\delta)\wedge\tau_r}(x)-x|\ge
		\lambda\right]\le 
		\lambda^{-2}(2\rho_0(r)^2\Delta^2+8\rho_1(r)^2\Delta).
\end{equation*}
Now note that $B_\lambda(x)$ is strictly included in $B_{|x|+2\lambda}$, 
so that the first exit time from $B_\lambda(x)$ is no later than
$\tau_{|x|+2\lambda}$.  But then the events
\begin{equation*}
	\left\{\omega:\sup_{0\le\delta\le\Delta}
	|\xi_{s,(s+\delta)\wedge\tau_{|x|+2\lambda}}(x)-x|\ge\lambda\right\},
	\quad
	\left\{\omega:\sup_{0\le\delta\le\Delta}
	|\xi_{s,s+\delta}(x)-x|\ge\lambda\right\}
\end{equation*}
are equivalent; after all, the events are equivalent on 
$\tau_{|x|+2\lambda}>s+\Delta$ by construction, whereas if
$\tau_{|x|+2\lambda}\le s+\Delta$ both events must be true as
$|\xi_{s,\tau_{|x|+2\lambda}}(x)-x|\ge\lambda$.  Hence
\begin{equation*}
	\mathbb{P}
	\left[\sup_{0\le\delta\le\Delta}
		|\xi_{s,s+\delta}(x)-x|\ge
		\lambda\right]\le 
		\lambda^{-2}(2\rho_0(|x|+2\lambda)^2\Delta^2
			+8\rho_1(|x|+2\lambda)^2\Delta)
\end{equation*}
where we have set $r=|x|+2\lambda$.  This completes the proof.
\qquad
\end{proof}

We now turn to the proof of the main theorem.

{\em Proof of Theorem \ref{thm:main}}.
	Let $\varepsilon>0$ and $Z=\{x\in\mathbb{R}^n:|x|>\varepsilon\}$.
	We begin by applying Lemma \ref{lemma:stochliou}.  To this
	end, define $S_\ell=\{x\in\mathbb{R}^n:|x|>\ell^{-1}\}$, so
	$S=\bigcup_\ell S_\ell=\mathbb{R}^n\backslash\{0\}$.  Clearly
	$D$ is integrable on $Z$ and there exists a 
	$C^2(\mathbb{R}^n,\mathbb{R}_+)$-approximation $D_\ell$ of $D$
	for each $\ell$.  It remains to check that $\tau_\ell\to-\infty$.
	Suppose that this is not the case; then given $x\in S$ there must 
	be a positive probability that $\xi_{s,t}^{-1}(x)=0$ for some 
	$-\infty<s<t$.  But $\xi_{s,t}^{-1}(0)=0$ for all $s$ and
	a.s.\ $\xi_{s,t}^{-1}(x)$ is one-to-one for all $s<t$, so this 
	cannot happen. Hence all the conditions of Lemma 
	\ref{lemma:stochliou} are satisfied, and we have
	\begin{equation}
	\label{eq:liouliou}
		0\le \int_ZD(x)\,dx+
		\int_s^t\mathbb{E}
		\int_{\xi_{\sigma,t}^{-1}(Z)}
			\mathscr{L}^*D(x)\,dx\,d\sigma.
	\end{equation}
	Now note that (\ref{eq:liouliou}) is nonincreasing with decreasing 
	$s$ due to $\mathscr{L}^*D\le 0$ and is finite because $D$ is 
	integrable on $Z$.  By monotone convergence the limit as 
	$s\to-\infty$ exists and is finite. Hence
	\begin{equation*}
		\int_{-\infty}^t\mathbb{D}(\xi_{\sigma,t}^{-1}(Z))\,d\sigma
		<\infty,
		\qquad\qquad
		\mathbb{D}(A)=
		-\int_{A}\mathscr{L}^*D(x)\,(\mathbb{P}(d\omega)\times
			\mu(dx))
	\end{equation*}
	where we have used Tonelli's theorem to convert the iterated
	integral to a single integral with respect to the product measure,
	and we slightly abuse our notation by writing
	$\mathbb{D}(\xi_{\sigma,t}^{-1}(Z))=\mathbb{D}(\{(\omega,x)\in
	\Omega\times S:\xi_{\sigma,t}(x,\omega)\in Z\})$.
	Note that $\mathscr{L}^*D\le 0$ implies that $\mathbb{D}$ is a 
	measure on $\Omega\times S$, and $\mathbb{D}$ is $\sigma$-finite 
	as $\mathbb{D}(\Omega\times\{x\in S:\tfrac{1}{k}<|x|<k\})<\infty$ 
	for all $k>1$ and $\bigcup_{n=2}^\infty
	(\Omega\times\{x\in S:\tfrac{1}{k}<|x|<k\})=\Omega\times S$.

	We now fix some $m\in\mathbb{N}$ and divide the halfline into bins 
	$S_k^m=[(k-1)2^{-m},k2^{-m}]$, $k\in\mathbb{N}$.  From each bin 
	we choose a time $t_k^m\in S_k^m$ such that
	\begin{equation*}
		\mathbb{D}(\xi_{t-t_k^m,t}^{-1}(Z))\le\inf_{s\in S_k^m}
			\mathbb{D}(\xi_{t-s,t}^{-1}(Z))+2^{-k}.
	\end{equation*}
	For fixed $m$, we denote this discrete grid by 
	$T_m=\{t_k^m:k\in\mathbb{N}\}$.  We now have
	\begin{equation*}
		2^{-m}\sum_{k=1}^\infty\mathbb{D}(\xi_{t-t_k^m,t}^{-1}(Z))
		\le 2^{-m}+
		\int_{-\infty}^t\mathbb{D}(\xi_{\sigma,t}^{-1}(Z))\,d\sigma
		<\infty.
	\end{equation*}
	Using the fact that the transformation $\theta_t$ of Theorem
	\ref{thm:theta} is $\mathbb{P}$-preserving to shift the times
	$t_k^m$ to the forward variable, we obtain
	\begin{equation*}
		\sum_{k=1}^\infty\mathbb{D}(\xi_{s,s+t_k^m}^{-1}(Z))
		=\sum_{k=1}^\infty\mathbb{D}(\xi_{t-t_k^m,t}^{-1}(Z))
		<\infty.
	\end{equation*}
	As $\mathbb{D}$ is $\sigma$-finite we can now apply the Borel-Cantelli
	lemma, which gives
	\begin{equation*}
		\mathbb{D}\left(\limsup_{k\to\infty}
			\xi_{s,s+t_k^m}^{-1}(Z)\right)
		=(\mathbb{P}\times\mu)\left(\limsup_{k\to\infty}
			\xi_{s,s+t_k^m}^{-1}(Z)\right)=0
	\end{equation*}
	where the first equality follows as $\mathscr{L}^*D(x)<0$ 
	$\mu$-a.e.\ implies $\mathbb{P}\times\mu\ll\mathbb{D}$.
	Consequently
	\begin{equation*}
		(\mathbb{P}\times\mu)\left(\bigcup_{m=1}^\infty\limsup_{t\in 
			T_m}\xi_{s,s+t}^{-1}(Z)\right)
		\le \sum_{m=1}^\infty
		(\mathbb{P}\times\mu)\left(\limsup_{t\in T_m}
			\xi_{s,s+t}^{-1}(Z)\right)=0.
	\end{equation*}
	We have thus shown that for all initial states $x$, except in
	a set $N\subset\mathbb{R}^n$ of Lebesgue measure zero, there is 
	$\mathbb{P}$-a.s.\ for any $m$ only a finite number of times $t$ 
	in the discrete grid $T_m$ such that $\xi_{s,s+t}(x)\in Z$.

	Let us fix an $x\not\in N$.  We now claim that the fact that 
	$\mathbb{P}$-a.s.\ for any $m$ there is only a finite number of 
	times $t\in T_m$ such that $\xi_{s,s+t}(x)\in Z$ implies that
	$\mathbb{P}$-a.s.\ $\limsup_{t\to\infty}|\xi_{s,t}(x)|\le\varepsilon$.
	To see this, suppose $\mathbb{P}[\limsup_{t\to\infty}
	|\xi_{s,t}(x)|>\varepsilon]=\delta>0$.  By monotone convergence
	$\mathbb{E}[\chi_{\limsup|\xi_{s,t}(x)|>\varepsilon'}]
	\nearrow\delta$ as $\varepsilon'\searrow\varepsilon$, hence there 
	exists an $\varepsilon'>\varepsilon$ such that
	$\mathbb{P}[\limsup_{t\to\infty}|\xi_{s,t}(x)|>\varepsilon']>0$.
	We have already shown, however, that a.s.\ $|\xi_{s,t}(x)|\le
	\varepsilon$ for infinitely many times $t_n\nearrow\infty$.
	Hence
	$$
	\mathbb{P}\left[\limsup_{t\to\infty}|\xi_{s,t}(x)|>\varepsilon'
		\right]>0
	\quad\Longrightarrow\quad
	\mathbb{P}[|\xi_{s,t}(x)|
		\mbox{ crosses }\varepsilon\mbox{ and }
		\varepsilon'\mbox{ infinitely often}]>0.
	$$
	Once we disprove latter statement, the claim is proved by 
	contradiction.

	To this end, introduce the following sequence of predictable 
	stopping times.  Let $\sigma_0=\inf\{t>s:|\xi_{s,t}(x)|
	\le\varepsilon\}$, $\tau_0=\inf\{t>\sigma_0:|\xi_{s,t}(x)|
	\ge\varepsilon'\}$, and for any $n>0$ we set
	$\sigma_n=\inf\{t>\tau_{n-1}:|\xi_{s,t}(x)|\le\varepsilon\}$,
	$\tau_n=\inf\{t>\sigma_n:|\xi_{s,t}(x)|\ge\varepsilon'\}$.
	Define 
	\begin{equation*}
		\Omega_n(\Delta)=\{\omega\in\Omega:\tau_n<\infty,~
			|\xi_{s,\tau_n+\delta}(x)|>\varepsilon~
			\forall\, 0\le\delta\le\Delta\}.
	\end{equation*}
	For any $\Delta>0$, the set of $\omega\in\Omega$ such that 
	$\omega\in\Omega_n(\Delta)$ for infinitely many $n$ must be of
	$\mathbb{P}$-measure zero; after all, we can choose $m$ 
	sufficiently large so that every time interval of length $\Delta$ 
	contains at least one point in $T_m$, and for points $t\in T_m$ we 
	have $|\xi_{s,t}(x)|>\varepsilon$ only finitely often 
	$\mathbb{P}$-a.s.  Thus $\sum_n\chi_{\Omega_n(\Delta)}<\infty$ 
	$\mathbb{P}$-a.s.
	To proceed, we use the following argument (see \cite{loeve}, pp.\ 
	398--399).  Introduce the discrete filtration 
	$\mathcal{B}_{k}=\mathcal{F}_s^{\tau_{k+1}}$ and define
	$Z_k=X_k-Y_k$ with
	\begin{equation*}
		X_k={\textstyle\sum}_{n=1}^k\chi_{\Omega_n(\Delta)},\qquad
		Y_k={\textstyle\sum}_{n=1}^k\mathbb{E}[\chi_{\Omega_n(\Delta)}|
			\mathcal{B}_{n-1}].
	\end{equation*}
	As $\Omega_k(\Delta)\in\mathcal{B}_{n}$ for all $k\le n$, $Z_k$ is 
	a $\mathcal{B}_k$-martingale. Now define for $a>0$ the stopping 
	time $\kappa(a)=\inf\{n:Z_n>a\}$.  As $|Z_k-Z_{k-1}|\le 1$ a.s., 
	the stopped process $Z'_k=Z_{k\wedge\kappa(a)}$ is a martingale 
	that is bounded from above, and by the martingale convergence 
	theorem $Z'_k$ converges a.s.\ as $k\to\infty$ to a finite random 
	variable $Z_\infty'$. But as $Z'_k$ and $Z_k$ coincide on 
	$\{\omega:\sup_nZ_n<a\}$ and $a>0$ was chosen arbitrarily, we conclude 
	that $Z_k\to Z_\infty<\infty$ on $\{\omega:\sup_nZ_n<\infty\}$ 
	(modulo a null set).  Note, however, that $X_n$ and $Y_n$ are both 
	positive increasing processes and we have already established that 
	$\sup_nX_n<\infty$ $\mathbb{P}$-a.s., so $\sup_nZ_n<\infty$ 
	$\mathbb{P}$-a.s.  But this implies that $Z_k$, and hence also 
	$Y_k$, converges to a finite value $\mathbb{P}$-a.s.  Thus we
	have established
	\begin{equation*}
		{\textstyle\sum}_{n=1}^\infty
		\mathbb{E}[\chi_{\Omega_n(\Delta)}|\mathcal{F}_s^{\tau_n}]
		<\infty\quad\mathbb{P}\mbox{-a.s.}\quad
		\mbox{for any }\Delta>0.
	\end{equation*}
	Note that by the continuity of the sample paths
	$|\xi_{s,\tau_n}(x)|=\varepsilon'$ on $\tau_n<\infty$.
	By Lemma \ref{lemma:uniform}, we can choose $\overline{\Delta}>0$ 
	sufficiently small such that
	\begin{equation*}
		P(y)=
		\mathbb{P}\left[\sup_{0\le\delta\le
			\overline{\Delta}}|\xi_{s,s+\delta}(y)-y|<
			\frac{\varepsilon'-\varepsilon}{2}
		\right]\ge\frac{1}{2}
	\end{equation*}
	for all $|y|=\varepsilon'$.   Using the strong Markov property,
	we can write
	\begin{equation*}
		\infty>
		{\textstyle\sum}_{n=1}^\infty
		\mathbb{E}[\chi_{\Omega_n(\overline{\Delta})}|
			\mathcal{F}_s^{\tau_n}]
		\ge
		{\textstyle\sum}_{n=1}^\infty
		P(\xi_{s,\tau_n}(x))\,\chi_{\tau_n<\infty}
		\ge
		\tfrac{1}{2}
		{\textstyle\sum}_{n=1}^\infty
		\chi_{\tau_n<\infty}~\mathbb{P}\mbox{-a.s.}
	\end{equation*}
	But this implies that $\tau_n<\infty$ finitely often 
	$\mathbb{P}$-a.s., contradicting the assertion that
	$\mathbb{P}[|\xi_{s,t}(x)|
	\mbox{ crosses }\varepsilon\mbox{ and } \varepsilon'
	\mbox{ infinitely often}]>0$.  This is the desired result.

	We have now shown that for $\mu$-almost all $x\in\mathbb{R}^n$, 
	$\mathbb{P}$-a.s.\ 
	$\limsup_{t\to\infty}|\xi_{s,t}(x)|\le\varepsilon$, i.e.\ 
	for $\mu$-almost all $x\in\mathbb{R}^n$ $\mathbb{P}$-a.s., 
	$\exists t_e>s$ such that $|\xi_{s,t}(x)|\le\varepsilon$ for $t\ge 
	t_e$.  But as this holds for any $\varepsilon>0$ the flow must 
	converge to the origin.  
\qquad\endproof

The proof of Theorem \ref{thm:main} is readily extended to prove other 
assertions, such as the following instability theorem.

\begin{theorem}
\label{thm:escape}
Let $X_k$, $k=0\ldots m$ be globally Lipschitz continuous and 
$C^{2,\alpha}$ for some $\alpha>0$.  Suppose there exists a $D\in 
C^2(\mathbb{R}^n,\mathbb{R}_+)$ such that $\mathscr{L}^*D(x)<0$ for 
$\mu$-almost all $x$.  Then for every initial time $s$ and $\mu$-almost 
every initial state $x$ the flow escapes to infinity, i.e.\ 
$|\xi_{s,t}(x)|\to\infty$ as $t\to\infty$ $\mathbb{P}$-a.s.
\end{theorem}

\begin{proof}
	Let $\varepsilon>0$ and $Z'=\{x\in\mathbb{R}^n:|x|<\varepsilon\}$.
	Again we begin by applying Lemma \ref{lemma:stochliou}.  We can
	simply choose $S_\ell=S=\mathbb{R}^n$ for all $\ell$; by 
	non-explosion $\tau_\ell=-\infty$ and the remaining conditions
	are evident.  Hence
	\begin{equation*}
		0\le \int_{Z'}D(x)\,dx+
		\int_s^t\mathbb{E}
		\int_{\xi_{\sigma,t}^{-1}(Z')}
			\mathscr{L}^*D(x)\,dx\,d\sigma.
	\end{equation*}
	Proceeding in exactly the same way as in the proof of Theorem 
	\ref{thm:main} we can now show that for $\mu$-almost all 
	$x\in\mathbb{R}^n$, $\mathbb{P}$-a.s.\ 
	$\liminf_{t\to\infty}|\xi_{s,t}(x)|\ge\varepsilon$, i.e.\ 
	for $\mu$-almost all $x\in\mathbb{R}^n$ $\mathbb{P}$-a.s., 
	$\exists t_e>s$ such that $|\xi_{s,t}(x)|\ge\varepsilon$ for $t\ge 
	t_e$.  But as this holds for any $\varepsilon>0$ the flow must 
	escape to infinity. \qquad
\end{proof}

{\em Remark.} At first sight the statements of Theorems \ref{thm:main} and
\ref{thm:escape} may seem contradictory, but this is not the case.  The
essential difference between the theorems is the region in $\mathbb{R}^n$
on which $D$ is integrable.  Roughly speaking, the idea behind the proofs
of Theorems \ref{thm:main} and \ref{thm:escape} is to show that if
$\mathscr{L}^*D<0$ $\mu$-a.e., then the solution of the It\^o equation can
only spend a finite amount of time in any region on which $D$ is
integrable. Hence in Theorem \ref{thm:main} the solution will attract to
the origin, whereas in Theorem \ref{thm:escape} the solution attracts to
infinity.

If we try to satisfy the conditions of Theorems \ref{thm:main} and 
\ref{thm:escape} simultaneously we will run into problems.  Suppose we 
have a nonnegative $D\in C^2(\mathbb{R}^n)$, as in Theorem 
\ref{thm:escape}, which is integrable as in Theorem \ref{thm:main}.
Then $D$ is a normalizable density function, i.e.\ we could normalize $D$ 
and interpret it as the density of the It\^o equation at some point in 
time.  But then $\mathscr{L}^*D<0$ would imply that the associated 
Fokker-Planck equation does not preserve normalization of the density.
Evidently Theorem \ref{thm:escape} can only be satisfied if $D$ is not 
integrable, whereas Theorem \ref{thm:main} requires $D$ to have a 
singularity at the origin.  See \S\ref{sec:examples} for examples.

\section{Further results}
\label{sec:further}

In this section we extend the main result to cases in which the global 
Lipschitz condition is not necessarily satisfied.  We first show that the 
result of Theorem \ref{thm:main} still holds if we can convert the 
coefficients of (\ref{eq:SDE}) to be globally Lipschitz continuous through 
a suitably chosen time transformation.  In particular, this allows us to 
treat the case that $X_k$, $k=0\ldots m$ and their first derivatives are 
polynomially bounded, provided that some additional integrability 
conditions on $D$ are satisfied.  We also extend the main result to 
the case in which the flow is restricted to an invariant subset of 
$\mathbb{R}^n$ with compact closure.

\begin{theorem}
\label{thm:main2}
Let $X_k:\mathbb{R}^n\to\mathbb{R}^n$, $k=0\ldots m$ be measurable and let
$X_k(0)=0$.  Suppose there is a strictly positive measurable map 
$c:\mathbb{R}^n\to (0,\infty)$ such that $c(x)$ and $c(x)^{-1}$ are 
locally bounded, and such that $c(x)X_0(x)$ and $\sqrt{c(x)}\,X_k(x)$, 
$k=1\ldots m$ are globally Lipschitz continuous and $C^{2,\alpha}$ for 
some $\alpha>0$. Suppose there exists $D:\mathbb{R}^n\backslash\{0\}\to 
\mathbb{R}_+$ such that $D(x)/c(x)$ is $C^2$, is integrable on 
$\{x\in\mathbb{R}^n:|x|>1\}$, and $\mathscr{L}^*D(x)<0$ for $\mu$-almost 
all $x$.  Then for every initial time $s$ and $\mu$-almost every initial 
state $x$ the solution $x_t$ of {\rm (\ref{eq:SDE})} tends to the origin 
as $t\to\infty$ $\mathbb{P}$-a.s.
\end{theorem}

\begin{proof}
	Consider the It\^o equation
	\begin{equation}
	\label{eq:cito}
		y_t = y_s+\int_s^t c(y_\tau)X_0(y_\tau)\,d\tau+
		\sum_{k=1}^m\int_s^t \sqrt{c(y_\tau)}\,X_k(y_\tau)\,dW_\tau^k.
	\end{equation}
	We will write $Y_0(y)=c(y)X_0(y)$, $Y_k(y)=\sqrt{c(y)}\,X_k(y)$ 
	($k=1\ldots m$), and $\tilde D(y)=D(y)/c(y)$.  Note that by 
	construction $\mathscr{\tilde L}^*\tilde D(y)=\mathscr{L}^*D(y)$,
	where $\mathscr{\tilde L}^*$ is the adjoint generator of 
	(\ref{eq:cito}).  By our assumptions we can apply Theorem 
	\ref{thm:main} to the It\^o equation (\ref{eq:cito}).  Thus for 
	all $y_s\in\mathbb{R}^n$, except in a set $N$ with $\mu(N)=0$, 
	$y_t\to 0$ as $t\to\infty$ $\mathbb{P}$-a.s.

	Now choose any $y_s\not\in N$ and define
	\begin{equation*}
		\beta_t=\int_s^t c(y_\tau)\,d\tau,
		~~~~~~~ ~~~~~~~ \alpha_t=\inf\{s:\beta_s>t\}.
	\end{equation*}
	Note that $\alpha_\tau$ is an $\mathcal{F}_s^t$-stopping time
	for each $\tau$.  We claim that $\beta_t\to\infty$ as $t\to\infty$ 
	a.s.; indeed $y_t$ a.s.\ spends an infinite amount of time in an 
	arbitrarily small neighborhood of the origin, and as $c(x)^{-1}$ 
	is locally bounded $c(x)\ge\delta>0$ in any such neighborhood.  
	Moreover, $\beta_t<\infty$ a.s.\ for any $t$ as $c(x)$ is locally
	bounded, and hence $\alpha_t\to\infty$ as $t\to\infty$ a.s.
	From \cite{rogersw}, \S{V.26} (pp.\ 175) it follows that the time 
	rescaled solution $y_{\alpha_t}$ is equivalent in law to the 
	solution $x_t$ of (\ref{eq:SDE}).  But as almost all paths of the 
	process $y_t$ go to zero asymptotically and as 
	$\alpha_t\to\infty$ a.s., the result follows.
	\qquad
\end{proof}

\begin{corollary}
\label{cor:main2}
	Suppose that all the conditions of Theorem \ref{thm:main}
	are satisfied except the global Lipschitz condition.
	Suppose that additionally $X_k$, $k=1\ldots m$ satisfy
	$|X_k(x)|\le C_k(1+|x|^{p+1})$, $|\partial X_k(x)/\partial x^i|
	\le C_k'(1+|x|^p)$, and $|X_0(x)|\le C_0(1+|x|^{2p+1})$,
	$|\partial X_0(x)/\partial x^i|\le C_0'(1+|x|^{2p})$ for
	some $p\ge 1$ and positive constants $C_k,C_k'<\infty$.
	If $(1+|x|^p)^2D(x)$ is integrable on $\{x\in\mathbb{R}^n:|x|>1\}$,
	then Theorem \ref{thm:main} still holds.
\end{corollary}

\begin{proof}
Let $c(x)=(1+|x|^p)^{-2}$, and note that $c(x)$ is smooth, strictly 
positive and that $c(x)$ and $c(x)^{-1}$ are locally bounded.
Let $Y_0(x)=c(x)X_0(x)$ and $Y_k(x)=\sqrt{c(x)}\,X_k(x)$, $k=1\ldots m$;
as $\sqrt{c(x)}$ is smooth and $X_k$, $k=0\ldots m$ are $C^{2,\alpha}$, 
the coefficients $Y_k$ are also $C^{2,\alpha}$.  To show that $Y_k$ are 
globally Lipschitz continuous it suffices to show that their first 
derivatives are bounded.  Let us write $y_{k,i}(x)=
|\partial Y_k/\partial x^i|(x)$.  First consider the case $k\ge 1$.  Then
$$
	y_{k,i}(x)\le \frac{1}{1+|x|^p}\,\left|
		\frac{\partial X_k(x)}{\partial x^i}
	\right|+\frac{p|x|^{p-2}|x^i|}{(1+|x|^p)^2}\,|X_k(x)|\le
	C_k'+C_k\frac{p|x|^{p-2}|x^i|}{(1+|x|^p)^2}(1+|x|^{p+1})
$$
which is bounded.  Similarly,
$$
	y_{0,i}(x)\le 
	C_0'\frac{1}{(1+|x|^p)^2}(1+|x|^{2p})
	+C_0\frac{2p|x|^{p-2}|x^i|}{(1+|x|^p)^3}
		(1+|x|^{2p+1})
$$
is bounded.  Finally, $D(x)/c(x)=(1+|x|^p)^2D(x)$ is $C^2$, as $D(x)$ is 
$C^2$ and $c(x)^{-1}$ is smooth on $\mathbb{R}^n\backslash\{0\}$.  Hence, 
provided $D(x)/c(x)$ is integrable on $\{x\in\mathbb{R}^n:|x|>1\}$, we can 
apply Theorem \ref{thm:main2}.
\qquad
\end{proof}

We now turn our attention to stochastic differential equations which 
evolve on an invariant set.  The following notion of invariance is 
sufficient for our purposes:

\begin{definition}
	A set $K$ is called backward invariant with respect to 
	the flow $\xi_{s,t}$ if $\xi_{s,t}^{-1}(K)\subset K$ a.s.\ for all
	$s<t$.
\end{definition}

We can now formulate the following result.

\begin{theorem}
\label{thm:invset}
	Suppose the It\^o equation {\rm(\ref{eq:SDE})} evolves on a 
	backward invariant open set $K$ with compact closure 
	$\overline{K}$.  Let $X_k$, $k=0\ldots m$ be $C^{2,\alpha}$ for 
	some $\alpha>0$ and let $X_k(0)=0\in\overline{K}$.  Suppose there 
	exists $D\in C^2(\overline{K}\backslash\{0\},\mathbb{R}_+)$ such 
	that $\mathscr{L}^*D(x)<0$ for $\mu$-almost all $x\in K$.  Then 
	for every initial time $s$ and $\mu$-almost every initial state 
	$x\in K$ the flow $\xi_{s,t}(x)$ tends to the origin as 
	$t\to\infty$ $\mathbb{P}$-a.s.
\end{theorem}

\begin{proof}
	We will assume without loss of generality that $X_k$, $k=0\ldots m$
	are globally Lipschitz continuous.  Indeed, we can smoothly modify 
	$X_k$ outside $K$ to have compact support without changing the 
	properties of the flow in $K$, and as $X_k$ are already locally
	Lipschitz continuous their modifications will be globally 
	Lipschitz continuous.

	Let $\varepsilon>0$ and $Z=\{x\in K:|x|>\varepsilon\}$.
	$D$ is integrable on $Z$, as $D$ is bounded on $Z$ and $Z$ has 
	compact closure.  Let $S_\ell$ be an increasing sequence of
	open sets whose closure is strictly contained in $S=K\backslash\{0\}$,
	such that $\bigcup_\ell S_\ell=S$.  Then there exists a 
	$C^2(\mathbb{R}^n,\mathbb{R}_+)$-approximation $D_\ell$ of $D$
	for each $\ell$, obtained by smoothly modifying $D$ outside 
	$S_\ell$ so that its support is contained in $\overline{K}$.  That 
	$\tau_\ell\to-\infty$ follows from backward invariance and from 
	the one-to-one property of the flow.  Hence all the conditions of 
	Lemma \ref{lemma:stochliou} are satisfied, and we have
	\begin{equation*}
		0\le\int_{Z} D(x)\,dx
		+\int_s^t\mathbb{E}
		\int_{\xi_{\sigma,t}^{-1}(Z)}
		\mathscr{L}^*D(x)\,dx\,d\sigma.
	\end{equation*}
	The remainder of the proof proceeds along the same lines as
	Theorem \ref{thm:main}.
	\qquad
\end{proof}

\section{Examples}
\label{sec:examples}

~

{\em Example 1}.
Consider the It\^o equation
\begin{equation}
\label{eq:ex1}
\begin{split}
	dx_t&=(x_t^2-2x_t-z_t^2)\,dt+x_t\,dW_t, \\
	dz_t&=2z_t(x_t-1)\,dt+z_t\,dW_t.
\end{split}
\end{equation}
Note that the line $z=0$ is invariant under the flow of (\ref{eq:ex1}), 
where the solution $(x_t,0)$ for an initial state $(x_0,0)$ is given by
\begin{equation*}
	dx_t=(x_t^2-2x_t)\,dt+x_t\,dW_t.
\end{equation*}
This equation has an explicit solution (see also \cite{arnold3} for a 
detailed analysis of the dynamical behavior of this system)
\begin{equation*}
	x_t=\frac{
		x_0 e^{-2t}e^{W_t-\frac{1}{2}t}
	}{1-x_0\int_0^t e^{-2s}e^{W_s-\frac{1}{2}s}ds}.
\end{equation*}
Clearly $x_t(\omega)$, $\omega\in\Omega$ explodes in finite time if
\begin{equation*}
	x_0>\left(
	\int_0^\infty e^{-2s}e^{W_s-\frac{1}{2}s}ds
	\right)^{-1}<\infty.
\end{equation*}
Hence the system (\ref{eq:ex1}) is certainly not globally stable.

Nonetheless, almost all points $(x_0,z_0)\in\mathbb{R}^2$ are attracted to 
the origin.  To show this, apply Corollary \ref{cor:main2} with
\begin{equation*}
	D(x,z)=\frac{1}{(x^2+z^2)^2},~~~~~~~ 
	~~~~~~~\mathscr{L}^*D(x,z)=-\frac{3}{(x^2+z^2)^2}<0.
\end{equation*}
Hence for a.e.\ $(x_0,z_0)\in\mathbb{R}^2$, a.s.\ $(x_t,z_t)\to(0,0)$ as 
$t\to\infty$.

{\em Example 2}.
Consider the It\^o equation
\begin{equation}
\label{eq:ex2}
\begin{split}
	dx_t&=12\,(2z_t-1)x_tz_t\,dt-\tfrac{1}{2}x_t\,dt+(1-2z_t)x_t\,dW_t, \\
	dy_t&=-\tfrac{1}{2}y_t\,dt+(1-2z_t)y_t\,dW_t, \\
	dz_t&=-12z_tx_t^2\,dt+2(1-z_t)z_t\,dW_t.
\end{split}
\end{equation}
Let $R_t=2z_t-2z_t^2-x_t^2-y_t^2$.  By It\^o's rule we have
\begin{equation*}
	dR_t=-4\,(1-z_t)z_tR_t\,dt+2\,(1-2z_t)R_t\,dW_t.
\end{equation*}
Evidently the ellipse $\{(x,y,z)\in\mathbb{R}^3:2z-2z^2-x^2-y^2=0\}$ is
invariant under (\ref{eq:ex2}).  Local uniqueness of the solution implies 
that the interior of the ellipse is also invariant.  Hence 
$K=\{(x,y,z)\in\mathbb{R}^3:2z-2z^2-x^2-y^2>0\}$ is a (backward) invariant 
set for the system (\ref{eq:ex2}).  Consider
\begin{equation*}
	D(x,y,z)=\frac{1}{z^2},~~~~~~~ 
	~~~~~~~\mathscr{L}^*D(x,y,z)=-\frac{12x^2}{z^2}<0~~\mu\mbox{-a.e.}
\end{equation*}
Hence by Theorem \ref{thm:invset} for a.e.\ $(x_0,y_0,z_0)\in K$, a.s.\ 
$(x_t,y_t,z_t)\to(0,0,0)$ as $t\to\infty$.  Note that $(0,0,0)$ is 
certainly not globally stable: it is easily verified that any point with 
$x_0=0$ and $z_0\ne 0$ is not attracted to $(0,0,0)$ a.s., as in this case 
$z_t$ has a constant nonzero mean.

{\em Example 3}.
We consider again the system (\ref{eq:ex2}), but now we are interested in 
the behavior of points in the invariant set 
$K'=\{(x,y,z)\in\mathbb{R}^3:2z-2z^2-x^2=0,~y=0\}$.  As $K'$ is not an 
open set we cannot apply Theorem \ref{thm:invset} to study this case.

Define the transformation $(x,z)\mapsto p=x/z$.
Note that $p$ is the stereographic projection of $(x,y,z)\in K'$ which 
maps $(0,0,0)\mapsto\infty$.  As the fixed point $(0,0,0)$ cannot be 
reached in finite time, we expect that the stereographic projection gives 
a well-defined dynamical system on $\mathbb{R}$.  Using It\^o's rule and 
$2z-2z^2-x^2=0$ we obtain the It\^o equation
\begin{equation*}
	dp_t=\left(\frac{3}{2}+\frac{20}{2+p^2}
		\right)p_t\,dt-p_t\,dW_t.
\end{equation*}
Note that this expression satisfies a global Lipschitz condition.  Now 
consider
\begin{equation*}
	D(p)=\sqrt{2+p^2},~~~~~~~ 
	~~~~~~~\mathscr{L}^*D(p)=-\frac{42}{(2+p^2)^{3/2}}<0.
\end{equation*}
Hence by Theorem \ref{thm:escape} for a.e.\ $p_0\in\mathbb{R}$, 
a.s.\ $p_t\to\infty$ as $t\to\infty$. This implies that the point
$(x,y,z)=(0,0,0)$ is almost globally stable in $K'$.

\section{Application to control synthesis}
\label{sec:control}

Consider an It\^o equation of the form
\begin{equation}
\label{eq:control}
	x_t = x+\int_s^t(X_0(x_\tau)+u_\tau \tilde X_0(x_\tau))\,d\tau+
		\sum_{k=1}^m\int_s^t X_k(x_\tau)\,dW_\tau^k
\end{equation}
where $\tilde X_0$ and $X_k$, $k=0\ldots m$ are $C^{2,\alpha}$ for some 
$\alpha>0$, $X_k(0)=0$ and $u_t$ is a scalar control input.  We consider 
instantaneous state feedback of the form $u_t=u(x_t)$ where $u(x)$ is 
$C^{2,\beta}$ for some $\beta>0$ and $u(0)=0$.  Then by Theorems 
\ref{thm:main}, \ref{thm:invset} or by Corollary \ref{cor:main2}, 
$x_t\to 0$ as $t\to\infty$ a.s.\ for a.e.\ $x_0$ if there exists a $D(x)$, 
with additional properties required by the appropriate theorem, such that
\begin{equation}
\label{eq:convexity}
\begin{split}
	\mathscr{L}^*D(x)=
	\frac{1}{2} \sum_{k=1}^m & \sum_{i,j=1}^n
		\frac{\partial^2}{\partial x^i\partial x^j}
		(X^i_k(x)X^j_k(x)D(x)) \\ &
	-\sum_{i=1}^n
		\frac{\partial}{\partial x^i}
		(X_0^i(x)D(x)+\tilde X_0^i(x)u(x)D(x))<0\qquad \mu\mbox{-a.e.}
\end{split}
\end{equation}
Note that this expression is affine in $D(x)$ and $u(x)D(x)$ and the set 
of functions $(D(x),u(x)D(x))$ which satisfy (\ref{eq:convexity}) is 
convex.  This fact has been used in the deterministic case by Prajna, 
Parrilo and Rantzer \cite{prajna1} to search for ``almost stabilizing'' 
controllers for systems with polynomial vector fields using convex 
optimization.  As the stochastic case enjoys the same convexity properties 
as the deterministic Theorem \ref{thm:determin} this approach can also be 
applied to find stabilizing controllers for stochastic nonlinear systems.  
Note that that convex optimization cannot be used to search for globally 
stabilizing controllers using LaSalle's theorem \cite{prajna1} as 
LaSalle's convergence criterion \cite{lasalle,deng2} is {\em not} convex.

The purpose of this section is to briefly outline the method of 
\cite{prajna1} for the synthesis of stabilizing controllers using convex 
optimization.  We will also discuss a simple example.

Suppose that $\tilde X_0$ and $X_k$, $k=0\ldots m$ are polynomial 
functions (the case of rational functions can be treated in a similar 
way.)  Consider $D(x)$ and $u(x)$ parametrized in the following way:
\begin{equation}
\label{eq:parametrized}
	D(x)=\frac{a(x)}{b(x)^\gamma},~~~~~~~ ~~~~~~~
	u(x)=\frac{c(x)}{a(x)}.
\end{equation}
Here $b(x)$ is a nonnegative polynomial which vanishes only at the origin,
$a(x)$ is a polynomial that is nonnegative in a neighborhood of the 
origin and is such that $u(x)$ is $C^{2,\beta}$, $c(x)$ is a polynomial 
that vanishes at the origin and $\gamma>0$ is a constant.  The orders of 
the polynomials and $\gamma$ can be chosen in such a way that $D(x)$ 
satisfies the integrability and growth requirements of Corollary 
\ref{cor:main2}.  For fixed $b(x)$ and $\gamma$ consider the expression
\begin{equation}
\label{eq:ineq}
	-b(x)^{\gamma+2}\mathscr{L}^*D(x)>0\qquad \mu\mbox{-a.e.}
\end{equation}
with $D(x)$ and $u(x)$ given by (\ref{eq:parametrized}) and 
$\mathscr{L}^*$ given by (\ref{eq:convexity}).  Then (\ref{eq:ineq}) 
is a polynomial inequality that is linear in the polynomial coefficients 
of $a(x)$ and $c(x)$.  Our goal is to formulate the search for these 
coefficients as a convex optimization problem.

Verifying whether (\ref{eq:ineq}) is satisfied comes down to testing 
nonnegativity of a polynomial (a nonnegative polynomial can only vanish on 
a finite set of points, and hence is positive $\mu$-a.e.)  This problem is 
known to be NP-hard in general; however, a powerful convex relaxation was 
suggested by Parrilo \cite{parrilo}.  Instead of testing (\ref{eq:ineq}) 
directly we may ask whether the polynomial can be written as a {\em sum of 
squares}, i.e.\ whether $-b(x)^{\gamma+2}\mathscr{L}^*D(x)=\sum_ip_i(x)^2$ 
for a set of polynomials $p_i(x)$.  The power of this relaxation comes 
from the fact that every sum of squares polynomial up to a specified order 
can be represented by a positive semidefinite matrix; hence the search for 
a sum of squares representation can be performed using semidefinite 
programming.  As (\ref{eq:ineq}) is convex in $a(x)$ and $c(x)$ the 
following is a convex optimization problem:
$$
	\mbox{Find polynomials }a(x)\mbox{, }c(x)\mbox{ such that }
	-b(x)^{\gamma+2}\mathscr{L}^*D(x)\mbox{ is sum of squares.}
$$
This type of problem, known as a sum of squares program, can be solved in 
a highly efficient manner using the software SOSTOOLS \cite{sostools}.
We refer to \cite{parrilo,sostools,prajna1} for further details on the 
computational technique.

{\em Remark.} Note that $a(0)$ and $c(0)$ depend only on the value of the 
constant coefficient of the polynomials $a(x)$ and $c(x)$.  Thus  
$c(0)=0$ can easily be enforced by fixing the constant coefficient of 
$c(x)$.  To make sure $a(x)$ is nonnegative near the origin and $u(x)$ 
does not blow up we can, for example, require $a(x)$ to be of the form 
$\lambda+p(x)$ with $\lambda>0$ and $p(x)$ is a sum of squares that 
vanishes at the origin.

Note that if the It\^o equation (\ref{eq:control}) evolves on an
invariant open set $K$ with compact closure then the sum of squares 
relaxation is overly restrictive.  A related relaxation that only 
guarantees polynomial nonnegativity on $K$ for the case that $K$ is a 
semialgebraic set is considered in e.g.\ \cite{prajna2}.

{\em Example.}
The following example is similar to an example in \cite{prajna1}.
Consider the It\^o equation
\begin{equation*}
\begin{split}
	dx_t&=(2x_t^3+x_t^2y_t-6x_ty_t^2+5y_t^3)\,dt+(x_t^2+y_t^2)\,dW_t, \\
	dy_t&=u_t\,dt-(x_t^2+y_t^2)\,dW_t.
\end{split}
\end{equation*}
We choose $b(x,y)=x^2+y^2$ and $\gamma=2.5$.  Using SOSTOOLS we find a 
solution with controller of order $3$ and a constant $a(x)$.  Note that 
these choices satisfy the integrability requirements of Corollary 
\ref{cor:main2}.  We obtain a stabilizing controller
\begin{equation*}
	u(x)=-2.7\,x^{3}+4.6\,x^{2}y-6.7\,xy^{2}-3.4\,y^{3}
\end{equation*}
where
\begin{equation*}
\begin{split}
	-(x^2+y^2)^{4.5}\mathscr{L}^*(&x^2+y^2)^{-2.5}=
	0.35\,y^{6}- 0.0015\,xy^{5} \\ &
	+ 0.6\,x^{2}y^{4}+
	0.0026\,x^{3}y^{3}+0.33\,x^{4}y^{2}+ 0.004\,x^{5}y+ 0.13\,x^{6}
\end{split}
\end{equation*}
is a sum of squares polynomial.

{\em Remark.}
A drawback of this method is that $b(x)$ and $\gamma$ must be fixed at 
the outset.  We have found that the method is very sensitive to the 
choice of $b(x)$ and $\gamma$ even in the deterministic case; often an 
unfortunate choice will cause the search to be infeasible.  Moreover it is 
not clear, even if there exists for polynomial $\tilde X_0$, $X_k$ a 
rational $u$ which almost globally stabilizes the system, that a rational 
$D$ can always be found that satisfies $\mathscr{L}^*D<0$.  Nonetheless 
the method can be successful in cases where other methods fail, and as 
such could be a useful addition to the stochastic nonlinear control 
engineer's toolbox.

\section*{Acknowledgments}
This work was performed in Hideo Mabuchi's group, and the author 
gratefully acknowleges his support.  The author would like to thank Anders 
Rantzer, Luc Bouten, Stephen Prajna, Paige Randall and especially Houman 
Owhadi for insightful discussions and comments.  The author is 
particularly thankful to an anonymous referee for pointing out a gap in the 
proofs and for his careful reading of the manuscript, which has 
significantly improved the presentation.

\bibliographystyle{siam}
\bibliography{Stabil}

\begin{thebibliography}{10}

\bibitem{arnold2}
{\sc L.~Arnold}, {\em The unfolding of dynamics in stochastic analysis}, Mat.
  Apl. Comput., 16 (1997), pp.~3--25.

\bibitem{arnold3}
\leavevmode\vrule height 2pt depth -1.6pt width 23pt, {\em Random Dynamical
  Systems}, Springer-Verlag, 1998.

\bibitem{baxendale}
{\sc P.~H. Baxendale}, {\em Properties of stochastic flows of diffeomorphisms},
  in Diffusion Processes and Related Problems in Analysis, Volume II, M.~A.
  Pinsky and V.~Wihstutz, eds., Progress in Probability 27, Birkh{\"a}user,
  1992.

\bibitem{bensoussan}
{\sc A.~Bensoussan}, {\em Stochastic Control of Partially Observable Systems},
  Cambridge University Press, 1992.

\bibitem{deng1}
{\sc H.~Deng and M.~Krsti{\'c}}, {\em Stochastic nonlinear stabilization---part
  {I}: a backstepping design}, Syst. Contr. Lett., 32 (1997), pp.~143--150.

\bibitem{deng2}
{\sc H.~Deng, M.~Krsti{\'c}, and R.~J. Williams}, {\em Stabilization of
  stochastic nonlinear systems driven by noise of unknown covariance}, IEEE
  Trans. Automat. Control, 46 (2001), pp.~1237--1253.

\bibitem{florch1}
{\sc P.~Florchinger}, {\em A universal formula for the stabilization of control
  stochastic differential equations}, Stoch. Anal. Appl., 11 (1993),
  pp.~155--162.

\bibitem{florch2}
\leavevmode\vrule height 2pt depth -1.6pt width 23pt, {\em Lyapunov-like
  techniques for stochastic stability}, SIAM J. Control Optim., 33 (1995),
  pp.~1151--1169.

\bibitem{florch4}
\leavevmode\vrule height 2pt depth -1.6pt width 23pt, {\em A stochastic
  {J}urdjevic-{Q}uinn theorem}, SIAM J. Control Optim., 41 (2002), pp.~83--88.

\bibitem{hasminskii}
{\sc R.~Z. Has'minski\u{\i}}, {\em Stochastic Stability of Differential
  Equations}, Sijthoff \& Noordhoff, 1980.

\bibitem{isidori}
{\sc A.~Isidori}, {\em Nonlinear Control Systems}, Springer-Verlag, third~ed.,
  1995.

\bibitem{kunita3}
{\sc H.~Kunita}, {\em On the decomposition of solutions of stochastic
  differential equations}, in Stochastic Integrals, D.~Williams, ed., Lecture
  Notes in Mathematics 851, Springer-Verlag, 1981.

\bibitem{kunita1}
\leavevmode\vrule height 2pt depth -1.6pt width 23pt, {\em Stochastic
  differential equations and stochastic flows of diffeomorphisms}, in {\'E}cole
  d'{\'E}t{\'e} de Probabilit{\'e}s de Saint-Flour XII, P.~L. Hennequin, ed.,
  Lecture Notes in Mathematics 1097, Springer-Verlag, 1984.

\bibitem{kunita2}
\leavevmode\vrule height 2pt depth -1.6pt width 23pt, {\em Stochastic Flows and
  Stochastic Differential Equations}, Cambridge University Press, 1990.

\bibitem{kushner1}
{\sc H.~J. Kushner}, {\em Stochastic Stability and Control}, Academic Press,
  1967.

\bibitem{lasalle}
\leavevmode\vrule height 2pt depth -1.6pt width 23pt, {\em Stochastic
  stability}, in Stability of Stochastic Dynamical Systems, R.~F. Curtain, ed.,
  Lecture Notes in Mathematics 294, Springer-Verlag, 1972.

\bibitem{loeve}
{\sc M.~Lo{\`e}ve}, {\em Probability Theory}, Van Nostrand, third~ed., 1963.

\bibitem{manzon}
{\sc P.~Monz{\'o}n}, {\em On necessary conditions for almost global stability},
  IEEE Trans. Automat. Control, 48 (2003), pp.~631--634.

\bibitem{prajna2}
{\sc A.~Papachristodoulou and S.~Prajna}, {\em On the construction of
  {L}yapunov functions using the sum of squares decomposition}, Proc. 41st IEEE
  CDC, 3 (2002), pp.~3482--3487.

\bibitem{parrilo}
{\sc P.~A. Parrilo}, {\em Semidefinite programming relaxations for
  semialgebraic problems}, Math. Prog. B, 96 (2003), pp.~293--320.

\bibitem{sostools}
{\sc S.~Prajna, A.~Papachristodoulou, P.~Seiler, and P.~A. Parrilo}, {\em
  {SOSTOOLS}: Sum of squares optimization toolbox for {MATLAB}. {U}ser's guide,
  version 2.00}.
\newblock Available from \texttt{http://www.cds.caltech.edu/sostools}, 2004.

\bibitem{prajna1}
{\sc S.~Prajna, P.~A. Parrilo, and A.~Rantzer}, {\em Nonlinear control
  synthesis by convex optimization}, IEEE Trans. Automat. Control, 49 (2004),
  pp.~310--314.

\bibitem{rantzer}
{\sc A.~Rantzer}, {\em A dual to {L}yapunov's stability theorem}, Syst. Contr.
  Lett., 42 (2001), pp.~161--168.

\bibitem{rogersw}
{\sc L.~C.~G. Rogers and D.~Williams}, {\em Diffusions, {M}arkov Processes and
  Martingales, Volume 2: {I}t{\^o} calculus}, Cambridge University Press,
  second~ed., 2000.

\end{thebibliography}

\end{document}